# POINCARÉ'S WORK ON CELESTIAL MECHANICS: PREDICTABILITY VERSUS DETERMINISM IN THE CONTEXT OF RESTRICTED THREE-BODY PROBLEM

Dr. **ROSÁRIO LAUREANO**

Instituto Universitário de Lisboa (ISCTE-IUL), ISTAR-IUL, Portugal

maria.laureano@iscte-iul.pt

Prof. Dr. **MANUEL ALBERTO M. FERREIRA**

Instituto Universitário de Lisboa (ISCTE-IUL), ISTAR-IUL, Portugal

manuel.ferreira@iscte-iul.pt

**ABSTRACT**

The publication of *Principia Mathematica* (1678) by Newton became known the celestial bodies' motion laws, which characterize the Classical Mechanics. Thereafter made sense to search about the movement of these bodies from known initial conditions, particularly in the *3*-body problem. In this paper it is exposed the influence of Poincaré's work (1880's) in this problem on the beginning of Deterministic Chaos Theory based on the development of several new qualitative tools. The application of such tools led him to the discovery of a special kind of behavior – the dynamical instability. Contrary to the common idea that this kind of behavior was first evidenced by Lorenz (1963), the Poincaré's theoretical research was sufficiently clear about the existence of chaotic deterministic behavior. Until the time of Poincaré, there was a tacit assumption that the uncertainty in the output does not arise from any randomness in the dynamical laws, since they are completely deterministic, but rather from the lack of the infinite accuracy in the initial conditions. In this paper it is emphasized that the issues of determinism and predictability are distinct.

**Key words:** dynamical instability, Chaos Theory, Determinism, unpredictability

1. INTRODUCTION

We belong to a world in which insignificant details can have a great impact. Very tiny changes in the starting conditions of an event can have substantial, sometimes even dramatic, effects on subsequent results. However, for centuries the prevailing view of the Universe was that it "runs like clockwork", and its running can be numerically predicted from a given set of initial conditions. As it is shown, this viewpoint is naive: it was verified in many real-world phenomena that small differences in the initial conditions of a process can have a significant effect on the final outcome. The first mathematical tools necessary to understand this kind of real-world phenomena were given by Jules Henri Poincaré (1854-1912) in his work on the restricted *3*-body problem in Celestial Mechanics which led to the field of chaotic dynamics in Mathematics.

As the science that studies the motion of celestial bodies in space, the Celestial Mechanics is undoubtedly one of the oldest sciences. Its golden age corresponds to the 18$^{th}$ and 19$^{th}$ centuries when the basic results of classical were obtained. From the Newton's work, both mathematicians and astronomers have worked in order to determine the movement of these bodies from known initial conditions, particularly in the *3*-body problem. The Earth-Moon-Sun system was the most important astronomical applications where none of the gravitational forces could be ignored. Besides, there were very practical and rather urgent needs in Astronomy and in Navigation to acquire accurate solutions of the *3*-body problem. Poincaré's work is a response to a challenge posed by King Oscar II of Sweden and Norway. However, although the Poincaré's theoretical research was sufficiently clear from the existence of chaotic deterministic behavior, the evidence to the scientific community provided by his work was only possible due to the use of a computer through the work of the mathematician and meteorologist Edward Norton Lorenz (1917-2008) almost 75 years later (1). In order to understand the importance of Poincaré's work in dynamical systems and the historical guidelines of the development of Deterministic Chaos Theory, it is essential to take into account the new tools that technology provided to the progress of Science. In turn, the concept of Deterministic Chaos Theory has extensive applications in current research in natural sciences, engineering, financial markets and information systems.

## 2. POINCARE'S CONTRIBUTION ON THE RESTRICTED *3*-BODY PROBLEM AND EMERGENCE OF DETERMINISTIC CHAOS THEORY

The *n*-body problem consists on determining the possible motions of *n* bodies moving under no influence other than that of their mutual gravitational attraction:

*"Given an initial set of data with the initial positions $s_i(0)$, masses $m_i$ and velocities $s_i´(0)$ of n bodies (i=1,2,…,n), with $s_i(0) \neq s_j(0)$ for all mutually distinct i and j, to determine the motions of the n bodies, and to find their positions at other times t, in accordance with the laws of classical mechanics."*

The first complete mathematical formulation of the *n*-body problem appeared at the Book 1 (*Philosophia Naturalis*) of *Principia Mathematica* by Isaac Newton (1642-1727), in which he discloses their laws of classical mechanics. Newton proved that celestial bodies can be modeled as mass points. Expressing the gravitational interactions which determine the motion of *n* mass points, and using his second law of motion, he obtained an initial-value problem of autonomous ordinary differential equations (ODE). The system of *n* mass points is an example of deterministic dynamical system (DDS). It is a continuous DDS defined by *n* second order ODE. It contains *6xn* variables, since each mass point is represented by *3* (Euclidean) space-components and *3* velocity-components.

Because it was very important to understand the movement of celestial bodies, many scientists have sought to determine whether the gravitational theory of Newton allowed a thorough knowledge of the movement of bodies based on their initial conditions. The study of the *n*-body problem by analytical methods, rather than the geometrical methods used by Newton, led to different research directions: to search for general theorems on the motion of *n* bodies or to search for approaches to a particular movement in a given period of time. Energy and angular momentum are conserved in a *n*-body system but it was hoped that other quantities, such as the integrals of the motion, might be conserved. Ernst Heinrich Bruns (1848-1919) showed in 1887 that there could be no conserved quantities which could be expressible as algebraic functions of the positions and velocities of the bodies (2-3). He proved that the 10 classical integrals/solutions (*3* for the positions of center of mass, *3* for the velocity of center of mass, *3* for the angular momentum and *one* for the energy) are the only algebraically linearly independent of this system with *6xn* degrees of freedom. This only allows the reduction to *6xn-10* variables.

For *n=2*, the number of variables is then reduced to *2*. In this case the DDS (of *2* ODE) is integrable. The *2*-body problem was completely solved by Johann Bernoulli in early 18[th] century: each body travels along a conic section which has a focus at the system's centre of mass. This problem is easy to solve but it is surprisingly difficult for *n≥3,* in spite of the reduction of system's number of variables given by the Bruns' work.

Historically, the case *n=3* is the most interesting and studied, mainly because the *3*-body problem models the Earth-Moon-Sun system. Bruns' work reduces the system's number of variables from *18* to *8*. However the difficulties on the 3-body problem, even in the restricted *3*-body problem where the mass of one of the bodies is infinitely small, can be demonstrated by the large number of papers on the subject since the mid-eighteenth century until the early twentieth century: more than 850 works, many of them of great mathematicians. Many of the early attempts to solve the *3*-body problem were quantitative in order to find (exact or approximate) explicit solutions for particular situations. Alexis-Claude Clairaut (1713-1765) was the first one to publish an approximate solution to the problem of the moon, in *"Du Système du Monde dans les Principe de la Gravitation Universelle"* (1747), using approximations to the solutions in the form of infinite series. Jean le Rond D'Alembert (1717-1783) and Clairaut developed in 1747 a long-standing rivalry, both attempted to analyze the problem in some degree of generality. However, D'Alembert stated clearly that he understood the classical mechanics based on metaphysical principles. He believed that Newton's laws of motion were logical necessities and forgot that Newton had established them from experimental evidence. Leonhard Euler (1707-1783) developed a method for calculate the perturbations in planetary motions, which was published in 1748 in *"Theoria et motuum Planetarum cometarum"*. Joseph-Louis Lagrange (1736-1813) showed in 1772 that the 3-body problem could be reduced to a DDS with seven variables (4). Carl Jacobi (1804-1851) found in 1843 another solution/integral to the problem (4). Together with William Rowan Hamilton (1795-1865), he was responsible for developing new methods of integration in ODE which were important in the later study carried out by Poincare. Charles-Eugène Delaunay (1816-1872) was the first one (1860) to provide a total elimination of secular terms in the problem of lunar theory (terms that grow infinitely or decrease with time and leading to an entirely new configuration of the system), but the series used to approximate the solutions converge very slowly (6). Another very important contribution is due to George William Hill (1838-1914), who published in 1877 and 1878 periodic solutions to the *3*-body problem (7-8). The knowledge of these periodic solutions greatly has influenced the subsequent Poincaré's research.

The question of solving the Newtonian n-body problem, with *n>2*, remained so difficult in the late 19th century that the King Oscar II, advised by Gösta Mittag-Leffler (the Majesty's science adviser), established in 1888 an award for its resolution. The award committee was constituted by the mathematicians Karl Weierstrass (1815-1897), Charles Hermite (1822-1901) and Mittag-Leffler (1846-1927). The first task of this committee was to define exactly the meaning of *"solving the n-body problem"*. Weierstrass carry out this task and the formulation, officially announced in *Acta Mathematica*, vol.VII, p.1885-1886, was the follows:

*"Given a system of arbitrarily many mass points that attract each other according to Newton's law, under the assumption that no 2 points ever collide, try to find a representation of the coordinates of each point as a series in a variable that is some known function of time and for all of whose values the series converges uniformly."*

Weierstrass's formulation was a reflection of the mainstream perception of his time. He equated the task of solving the *n*-body problem to the task of finding a power series solution that converges for all time. At this time, mathematicians had tried to solve ODE by deriving analytic formula for solutions. Since solutions in closed form are in general not attainable, power series were used as a substitute. Nowadays, modern computer can do the job required through numerical integration; but at the time people had to rely on power series solution and it was not uncommon for a mathematician to spend his entire career computing power series solution of the *3*-body problem, using all the tricks one could possibly imagine in attaining more accurate solution that can better predict the position of a body. However, power series only converge on finite time interval, and the size of the convergence interval varies depending on the location of the solution. Therefore, it became questionable, since the latter half of the 19th century, if the traditional methods trying to solve the problem by series were sufficient.

Poincaré was one of the scientists who started working on the proposed challenge. Poincaré began his research into this subject with a completely different approach: the overall understanding of the behavior of all solutions/integrals of a DDS was more important than the approximate knowledge of the local behavior of solutions defined analytically. There he used his doctoral dissertation on Mathematics entitled *"Sur les propriétés des fonctions définies par les équations différences"* (1879, advised by Hermite) in which he devised a new way of studying the general geometric properties of functions defined by ODE, based on topology and using the Lobachevsky's non-Euclidean geometry. In the discussion of his dissertation he clearly said that it could be used to model the behavior of multiple bodies in free motion within the solar system, where the essential question is about the stability and qualitative properties of planetary orbits and not about the numerical solution. Poincaré's geometric point of view, which is by the way simple, is to view an n-vector as a point in the called phase space and to view solutions as a collection of non-intersecting curves in this space. He pointed out that, as a mathematical problem, the questions we should ask in the study of ODE are (i) what kind of solution curves are allowed by the ODE, and (ii) in what way all these solution curves fit together.

His brilliant insight in the resolution of the problem was to recognize that rather than considering the flow in the entire *3*-dimensional space, it was much more convenient to consider the called *first return map* (or *Poincaré map*) induced by the flow on a *2*-dimensional surface of section *S* transverse to the flow.

This map is defined by choosing a point *P* on *S* at which *S* is intersected by a flow line, the image of *P* under the map is then the point *P'* at which that flow line first intersects *S* again. If x(t, x(0)) is the solution of the ODE satisfying x(0, x(0)) = x(0), the first return map iterated is defined by x(0) →x(T, x(0)) in the 2-dimensional space. This map is the time-T map. To study the solutions of an ODE is the same as to study the iterations of a map that is defined by the solutions of this ODE. A periodic solution corresponds in the 3-dimensional space to a closed curve, but under the map a *2π-periodic* solution corresponds to a fixed point and a *2pπ-periodic* solution corresponds to a cycle of period *p*.

So, Poincaré turned his attention to the unstable periodic solutions and to the study of behavior of others solutions in the neighborhood of them. He found that there are *asymptotic solutions*, that is, solutions which slowly either approach or move away (asymptotically approach in both the backward and the forward time) from an unstable periodic solution. Poincaré called these asymptotic solutions as *homoclinic solutions*. He showed that in the *3*-dimensional solution space of the restricted *3*-body problem, these homoclinic solutions generate families of curves which fill out *invariant surfaces* in the phase space and which asymptotically approach the curve representing the generating periodic solution. These asymptotic surfaces correspond to curves in the transverse section and, in order to understand the behavior of asymptotic solutions, he investigated the nature of these curves on transverse section.

When Poincaré sought to highlight the properties of the exact equations of asymptotic surfaces, he constructs the intersection curves *AO'B'* and *A'O'B* of these surfaces with a transversal section. The unstable periodic solution was represented by a closed trajectory which intersects the transverse section at the point *O'*. Poincaré considered that the points *B* and *B'* coincide, then the asymptotic surfaces were closed. So, the curve *BO'B'* delimits a closed surface. Therefore, an unstable periodic solution corresponded to a system of asymptotical solutions confined to a region of space. This conclusion was taken by Poincaré as a result of stability. Poincaré had used the definition of stability known as Poisson stability. In a letter addressed to Mittag-Leffler, he claimed that he has proved a stability result for the restricted 3-body problem. He wrote (9):

*"In this particular case, I have found a rigorous proof of stability and a method of placing precise limits on the elements of the third body… I now hope that I will be able to attack the general case and … if not completely resolve the problem (of this I have little hope), then at least found sufficiently complete results to send into the competition."*

The work was then submitted on May 17, 1888, and the award was presented to him at January of 1889. In total, 12 papers were sent anonymously to the competition's judges. The Poincaré's work was dually refereed by Mittag-Leffler and Weierstrass, and the latter asserted in his report to the former that (9):

*"I have no difficulty in declaring that the memoir in question deserves the prize. You may tell your Sovereign that this work cannot, in truth, be considered as supplying a complete solution to the question we proposed, but it is nevertheless of such importance that its publication will open a new era in the history of celestial mechanics. His Majesty's goal in opening the contest can therefore be considered attained."*

However, in the middle of trying to answer questions raised by Edvard Phragmén (1863-1937) (editor of *Acta Mathematica* and responsible for copyediting his paper for publication), Poincaré soon realized that his work contained a mathematical error and the stability result he claimed to Mittag-Leffler was wrong (July of 1889). But a limited number of copies of journal with the mistaken manuscript were already in print and mailed out. Mittag-Leffler decided to make a recall and he ordered all recalled copies physically destroyed (note that Mittag-Leffler founded the journal *Acta Mathematica* in 1882, precisely with the help of King Oscar's sponsorship). Poincaré spent his entire prize money and some to pay for the double printing.

While trying to detect the mathematical consequences of his mistake, Poincaré proved that the problem did not have a solution because the existence of the called *homoclinic tangle*. He verifies that it was wrong to consider the coincidence of the points *B* and *B'*. When he applies his corrected version to the asymptotic surfaces's study, he concluded that, for each unstable periodic trajectory, if the curves *AO'B'* and *A'O'B* have more than one common point, then there are an infinite number of intersections. Through these points of intersection are an infinite number of curves which, as belonging to both asymptotic surfaces, are called *doubly asymptotic trajectories*.

The trajectory of an homoclinic solution in phase space is an invariant loop for the time-T map of the unperturbed equation. Small perturbation, unfortunately, would usually break the loop into two intersecting curves (*homoclinic intersection*). What Poincaré discovered was when the loop is broken into two intersecting curves, which he knew by then as inevitable, then they would be forced to intersect each other to form a web – called the homoclinic tangle-, and the structure of this web appeared to be incomprehensibly complicated. For solutions in this complicated web, there would be no well-defined destination, therefore no possibility for dynamical stability.

In this way, Poincaré found the dynamical instability, currently designed by *chaotic behavior*. It is a special kind of behavior characterized by the existence of orbits in the system which are nonperiodic, and yet not forever increasing (exploding) to infinity nor approaching or orbiting around an equilibrium point, which indicated the lack of asymptotical stability in the system. The evolution of such a system is often chaotic in the sense that

*"It may happen that small differences in the initial conditions produce very great ones in the final phenomena. A small error in the former will produce an enormous error in the latter. Prediction becomes impossible.",*

as written by Poincaré in *Science and Method* **(10)** (1903). So, this finding gave birth to the Deterministic Chaos Theory, a mathematical theory of great influence in later times. Finally, Poincaré's research culminated in a paper published by Acta Mathematica (11) after the recall. The history has proved that it equaled all that was asserted in Weierstrass's report. Yet, it is made famous by the new tools of qualitative nature introduced. The discovery of such qualitative tools provided the future development of the Theory of Dynamical Systems.

The Deterministic Chaos Theory deals with the long-term (asymptotic) qualitative behavior of DDS. The focus is not on finding precise solutions to the equations defining the systems, which is often hopeless, but rather to answer questions like *(Q1)* Will the system settle down to a steady state in the long-term? *(Q2)* If so, what are the possible steady states of the system? *(Q3)* Does the long-term behavior of the system depends on its initial conditions?

These questions point out to the three most known properties of chaotic behavior: for a DDS to be classified as chaotic, it must have an invariant subset *E* of its phase space with the properties of *topological transitivity*, *density of periodic orbits* and *sensitive dependence on the initial conditions*.

### 3. DETERMINISM AND (UN)PREDICTABILITY

A common belief toward the end of the 19$^{th}$ century was that Mathematics, when properly combined with logic, can be used to obtain an exact description of the world around us. So, if we can know the state of the Universe at one moment, and write out all the laws that govern it, we can accurately predict its state at any other moment. This doctrine – the *Determinism* - states that the cause-and-effect rules completely govern all the motion and structure on the material level. So, in principle, every event can be completely predicted in advance, or in retrospect.

Determinism was inspired by Newtonian physics composed by a concise set of principles. Newton's laws, being expressed as mathematical equations and not simply as ordinary sentences, are completely deterministic. ODE of a given DDS connect the measurements at a given initial moment to the values at a later or earlier moment and allow to determine exactly what to expect, given a set of initial conditions. Two nearly-indistinguishable sets of initial conditions will always produce identically the same behavior at any moment in the future or past. Until the time of Poincaré, there was an implicit assumption among almost all scientists: if you could decrease the uncertainty in the real measurements then any imprecision in the prediction would decrease in the same way. This idea is called the "shrink-shrink" rule. Uncertainty in the final predictions was a minor problem because it was assumed that it does not arise from any properties of the dynamical laws since they are completely deterministic. It was theoretically possible to obtain nearly-perfect predictions for the behavior of any DDS: by putting more and more precise inputs, we got more precise outputs for any later or earlier moment.

Although certain DDS did indeed obey the "shrink-shrink" rule, the 3-body system did not. Poincaré's work was a proof that for some systems, like the DDS which models the restricted 3-body problem, the only way to obtain predictions with significant degree of accuracy is to determine the initial conditions with absolutely infinite precision, because one of the properties of these systems is the sensitivity to initial conditions. Any imprecision at all in the initial conditions, no matter how small, result after a short period of time in an uncertain deterministic prediction. If a chaotic output is generated by one set of initial conditions and then they are changed, even just a little bit like a perturbation, the output will change over time. This expansion (blowing up) of tiny uncertainties in the initial conditions into enormous uncertainties in the outcomes remained even if the initial uncertainties were decreased to smallest imaginable size. So, there was a conflict with the paradigm of Modern Science. That is, although the system follows deterministic rules, its time evolution appears random. The system is predictable in principle but unpredictable in practice. Even using perfect measuring devices, it is impossible to record a real measurement with infinite precision, since this requires displaying an infinite number of digits. In dynamics, the presence of uncertainty in any real measurement means that the initial conditions in studying any DDS cannot be specified to infinite accuracy.

## 4. CONCLUSIONS

In common language, chaos means disorder or randomness. In mathematics there is much more: it is a very specific kind of unpredictability, that is a deterministic behavior where the "shrink-shrink" rule is not valid. It is also common to think that chaotic DDS are random. However, they are deterministic systems governed by mathematical equations, and so they are completely predictable given perfect knowledge of the initial conditions. There is no inherent randomness in this regard. Such accurate knowledge is never attainable in real life because slight errors are intrinsic to any real measurement. The issues of Determinism and predictability are then distinct.

The discovery of Chaos had important consequences for modeling real-world systems. Given a DS exhibiting random behavior, a useful question is whether the DS is better approximated using a deterministic model that allows chaotic dynamics or, alternatively, by a stochastic model, or yet by a mixed deterministic/stochastic model. Chaos is the rule and not the exception. There are chaotic examples in economics, fluid dynamics, optics, chemistry, climate changes, biology, population dynamics, engineering, etc (14-19). There are also examples where Chaos can be used and/or manipulated for certain purposes. It is the case of chaotic encryption which is used to protect digital information, dynamical chaos control and chaos synchronization (20-21).